\documentclass[12pt]{amsart}
\usepackage{tikz}
\newcommand{\comments}[1]{}
\usepackage{amsmath, amssymb, amsthm, color, verbatim}
\usepackage{hyperref}

\newtheorem{theorem}{Theorem}[section]
\newtheorem{lemma}[theorem]{Lemma}

\theoremstyle{definition}
\newtheorem{example}[theorem]{Example}

\DeclareMathOperator{\CT}{\mathrm{{\textbf{CT}}}}

\DeclareMathOperator{\Char}{\mathrm{{\textbf{Char}}}}

\begin{document}

\title[Character Table Database]
{A Database of Groups with Equivalent Character Tables}
\author{William Cocke}
\address{William Cocke\\Department of Mathematics, University of Wisconsin-Madison}
\email{cocke@math.wisc.edu}
\author{Steve Goldstein}
\address{Steve Goldstein\\Department of Botany and Department of Mathematics, University of Wisconsin-Madison}
\email{sgoldstein@wisc.edu}
\author{Michael Stemper}
\address{Michael Stemper\\Department of Mathematics, University of Wisconsin-Madison}
\email{mstemper2@wisc.edu}

\maketitle

\begin{abstract}
Two groups are said to have the same character table if a permutation of the rows and a permutation of the columns of one table produces the other table. The problem of determining when two groups have the same character table is computationally intriguing. We have constructed a database containing for all finite groups of order less than 2000 (excluding those of order 1024), a partitioning of groups into classes having the same character table. To handle the 408,641,062 groups of order 1536 and other orders with a large number of groups we utilized high-throughput computing together with a new algorithmic approach to the problem. Our approach involved using graph isomorphism software to construct canoncial graphs that correspond to the character table of a group and then hashing the graphs. 
\end{abstract}

\section{Introduction} 
There are many questions about what properties of a group are captured by the character table. In general, the character table is like a shadow of the corresponding group: some information is preserved and other information is lost. Two groups may have the same character table, e.g., the dihedral and quaternion groups of order 8, or the two extra-special groups of order $p^3$ for any prime $p$. It is also easy to show that certain groups cannot have the same character table, e.g., groups of different orders or with a different number of conjugacy classes.

To aid in the study of what information about a group is preserved and what is
lost in the passage to its character table we have constructed a
database containing for all finite groups of order less than 2000
(excluding those of order 1024), a partitioning of groups into classes
having the same character table.

The database reveals
new computational examples of character tables not preserving certain properties of the underlying groups. 
For example, previous work by Mattarei gave examples of groups $G$ and $H$ with the same character table but different derived lengths \cite{Mattarei1,Mattarei2}. The smallest examples given by Mattarei had order $5^{11}$. Using our database we have the following computational observations.

\begin{example}
Let $G$ and $H$ be two groups with $|G|=|H|<512$. Suppose that the character tables of $G$ and $H$ are the same. Then the derived length of $G$ and $H$ are the same. \end{example}
\begin{example}
There are groups of order 512 that share a character table, but have different derived lengths. See Table 3 in Section \ref{sec derived length}.
\end{example}

The structure of the database itself is of some interest. For example,
we have expanded the table found in the book by Lux and Pahlings
\cite[Table 2.2 pg 136]{LP}.
See Tables 1 and 2.

\begin{table}[h!]
\begin{center}
\begin{tabular} {r  r r r  }
Order & Number of Groups & Number of Tables  & Largest Class Size\\
\hline
2 & 1 & 1 &1  \\
4 & 2 & 2 & 1 \\
8 & 5 & 4 & 2 \\
16 & 14 & 11 & 2  \\
32 & 51 & 35 & 3\\
64 & 267 & 146 & 6\\
128 & 2328 & 904 & 36\\
256 & 56092 & 9501 & 256\\
512 &  10494213  & 360134& 135424  
\end{tabular}
\caption{The number of groups and character tables of order $2^k$. The largest class size is the largest number of groups of that order that share a character table.}
\end{center}
\end{table}

\begin{table}[h!]
\begin{center}
\begin{tabular} {r  r r r  }
Order & Number of Groups & Number of Tables  & Largest Class Size\\
\hline
6 & 2 & 2 &1  \\
12 & 5 & 5 & 1 \\
24 & 15 & 13 &2 \\
48 & 52& 42 & 2  \\
96 & 231 & 160 & 3\\
192 & 1543 & 834 & 9\\
384  & 21185 & 7237 & 36\\
768 & 1090235 & 139974 & 512\\
1536&  408641062  & 20540010 & 135424  
\end{tabular}
\caption{The number of groups and character tables of order $2^k\cdot 3$. The largest class size is the largest number of groups of that order that share a character table.}
\end{center}
\end{table}

Table 1
was based on work of Skrzipczyk \cite{Skrzipczyk} in which she
searched for a minimal example of a Brauer pair. Skrzipczyk's approach
was to minimize the number of
pairwise
comparisons by clustering character tables
using 
group invariants that are known to be captured by the character
table. Skrzipczyk was able to partition the $2$-groups of order up to
256 by their character tables. However, this approach does not scale
well
because its complexity is quadratic in the size of the largest cluster. We will compare this approach with ours in Section \ref{sec compare}.

To handle the 408641062 groups of order 1536 and other orders with a large number of groups we utilize a novel approach to the problem. 
We clustered the groups of a given order by number of conjugacy 
classes. For each group within a cluster, we calculated a
string that faithfully represents the character table using an efficient 
implementation of a canonical graph labeling algorithm, effectively 
linearizing the pairwise comparison problem. We then verified there were 
no hash collisions. Hence, this 
hash-of-string-representation-of-the-canonically-labeled-graph is, 
empirically, an invariant of the group that determines its character table, or 
in the language of Section \ref{sec history}, its equitabular class.

Even though this hashing scheme enabled an efficient implementation for 
constructing the database, deploying it for hundreds of 
millions of groups was, nonetheless, a formidable challenge.  
We addressed this challenge by performing the calculations on a high-throughput computing pool. Indeed our hashing scheme was conceived with this resource in mind. 
For any given group, 
calculating the hash, the computationally-intensive step in the algorithm, can 
be a stand-alone computation, each one producing a small text file as 
output and having small processor and memory footprints, making it a 
good match for high-throughput computing.

The rest of this paper proceeds as follows. Section \ref{sec history} briefly recalls the definition of a character table, and what it means for two character tables to be equal. Section \ref{sec methods} contains the techniques used to produce the database. Section \ref{sec observations} includes some observations and applications from the database. Section \ref{sec availability} provides information on how to obtain the database and Section \ref{sec acknowledgements} contains a few acknowledgements relevant to the project. 

\section{The character table of a group.}\label{sec history}

For a finite group $G$, a character of $G$ is the trace of a representation $\phi:G\rightarrow \text{GL}_n(\mathbb{F})$, where $\text{GL}_n(\mathbb{F})$ is the general linear group of degree $n$ over the field $\mathbb{F}$. For general facts about character theory, including more definitions, we refer the reader to Isaacs \cite{Isaacs_CT}. In this paper we are only interested in characters over the complex numbers and will assume that $\mathbb{F}=\mathbb{C}$. A character is called irreducible if the corresponding representation is irreducible. The set of irreducible characters of a group $G$ is written as $\text{Irr}(G)$. We write $\textbf{Classes}(G)$ for the set of conjugacy classes of $G$. 
 
To define a character table of a group $G$, we must first fix orderings of $\text{Irr}(G)$ and $\mathbf{Classes}(G)$. A character table $\CT(G)$ of a group $G$ is an array whose $ij$-entry is the value of the $i$-th character in $\text{Irr}(G)$ evaluated on a representative of the $j$-th conjugacy class in $\mathbf{Classes}(G)$. Clearly, the presentation of $\CT(G)$ as an array depends on the orderings we choose for $\text{Irr}(G)$ and $\mathbf{Classes}(G)$. In general, we will say that two character tables are equivalent if there is a permutation of the rows and a permutation of the columns of one table that equals the other table and we write this as $\CT(G) \cong \CT(H)$. If $\CT(G)\cong \CT(H)$, we say that $G$ and $H$ are equitabular; this terminology is non-standard, but useful for the work in this paper. For a group $G$ we will call the set of all groups that are equitabular with $G$ the equiatabular class of $G$. Hereafter when we say two groups have the same character table, we mean that the character tables are equivalent in this manner.  

\sloppy Note that the same group might produce a large number of different character tables that are all equivalent. We can view the character table $\CT(G)$ of $G$ as a mapping 
\[\CT(G): \text{Irr}(G) \times \textbf{Classes}(G) \rightarrow \mathbb{C} \] where for a character $\chi$ and a conjugacy class $\mathcal{C}=g^G$, we have \[\CT(G)(\chi,\mathcal{C}) = \chi(g).\] Then two character tables $\CT(G)$ and $\CT(H)$ are equivalent if and only if there are bijections $\Phi:\text{Irr}(G)\rightarrow \text{Irr}(H)$ and $\Psi:\textbf{Classes}(G) \rightarrow \textbf{Classes}(H)$ such that \begin{equation} \label{eqn 1} CT(G)(\chi,\mathcal{C})=\CT(H)(\Phi(\chi),\Psi(\mathcal{C})).\end{equation} In Section \ref{sec step 2} we will see that equation \ref{eqn 1} is reflected in the graphs we construct.

\section{Constructing the database}\label{sec methods}

In this section we do two things: we explain our technique to construct the database and compare this technique against the standard approach with an emphasis on our reliance on high-throughput computing. In Sections \ref{sec step 1} to \ref{sec collisions} we discuss our algorithms to build the partition of groups by their equitabular classes. Sections \ref{sec step 1} and \ref{sec step 2} also give an example of how our encoding of $\CT(G)$ as a graph would work for symmetric group on 3 symbols. In Section \ref{sec compare} we compare our algorithm with previous approaches to the problem. 

To construct the equitabular class of a group $G$ we need only to examine groups with the same order and the same number of conjugacy classes. For two finite groups $G$ and $H$ the software package {\sc Magma}\cite{Magma} can be easily used to determine if $G$ and $H$ are equitabular. In many cases, we were able to partition the groups into equitabular classes by running one-versus-one comparisions using the order of the group and the number of conjugacy classes to separate the groups from one another. However for orders with a large number of groups, such well-chosen one-to-one comparison was not feasible.  We used the {\sc GAP} \cite{GAP} package {\sc Digraphs} \cite{Digraphs} along with { \sc HTCondor} \cite{Condor} to build the partition for orders with a larger number of groups. 

Our approach to handle orders with a large number of groups was as follows: 

\begin{itemize}
\item[(1)] Calculate the character table for a group $G$.
\item[(2)] Represent $\CT(G)$ as a graph $\mathcal{G}$.
\item[(3)] Find a canonical labeling of $\mathcal{G}$. 
\item[(4)] Hash the canonical labeling of $\mathcal{G}$. 
\item[(5)] Sort the list of hashes.
\item[(6)] Check for hash collisions.  
\end{itemize}

The following subsections cover each of the above steps. 
\subsection{Calculating a character table for a group $G$.}\label{sec step 1}

We used the native character table construction function within {\sc GAP} to construct character tables. For a group $G$, the function \texttt{CharacterTable(G)} returns $\CT(G)$.

Let $S$ be the symmetric group on 3 symbols. Then \[\CT(S)=\begin{pmatrix}
1 & 1 &1 \\
1 & -1 & 1\\
2 & 0 & -1
\end{pmatrix}.\]

\subsection{Represent $\CT(G)$ as a graph.}\label{sec step 2}

We want to represent the character table $\CT(G)$ as a graph $\mathcal{G}$. We will actually build a colored digraph $\mathcal{G}$ from the table representing $\CT(G)$; here the use of color merely serves to identify that two vertices cannot be interchanged by the graph. Note that coloring the vertices is equivalent to constructing a slightly more complicated graph that prevents the existence of automorphisms that interchange certain vertices. 

Suppose that $\CT(G)$ is an $n$-by-$n$ table with $k$ distinct entries. Then our graph $\mathcal{G}$ will have vertex set 
\[
V={\color{red}\{r_1,\dots,r_n\}} \cup {\color{blue}\{c_1,\dots,c_n\}} \cup {\color{cyan}\{e_{i,j}: 1\le i,j \le n \}} \cup \{v_1,\dots,v_k\},
\] where as indicated, the $r$-vertices are all one color, the $c$-vertices are a different color, and the $e$-vertices are a third color;  moreover each of the $v$-vertices is a distinct color, meaning that the entire vertex set $V$  is colored by $3+k$ colors. The $r$-vertices will correspond to rows of $\CT(G)$, the $c$-vertices to columns, the $e$-vertices to entries, and the $v$-vertices to distinct values of $\CT(G)$. The edge set of $\mathcal{G}$ is 
\begin{flalign*}
E=&\left(\bigcup_{1\le i \le n} \{(r_i,e_{i,j}) : 1\le j \le n\} \right) \cup \left(\bigcup_{1\le j \le n} \{(c_j,e_{i,j}): 1\le i \le n\} \right) \\
&\cup \left(\bigcup_{1\le \ell \le k} \{(v_\ell, e_{i,j}): \CT(G)[i,j] = v_\ell, 1\le i,j \le n\} \right).\\
\end{flalign*}

Using $\CT(S)$ as written in Section \ref{sec step 1}, the graph $\mathcal{S}$ is below. 

\begin{figure}
\begin{center}
\begin{tikzpicture}[scale=.5]
\node[red, draw,circle,inner sep=2pt, fill] at (4,3) {};
\node[red, draw,circle,inner sep=2pt, fill] at (4,5) {}; 
\node[red, draw,circle,inner sep=2pt, fill] at (4,7) {};
\node[draw,circle,inner sep=1.5pt, fill] at (7,2) {};
\node[draw,circle,inner sep=1.5pt, fill] at (7,4) {};
\node[draw,circle,inner sep=1.5pt, fill] at (7,6) {};
\node[draw,circle,inner sep=1.5pt, fill] at (9,2) {};
\node[draw,circle,inner sep=1.5pt, fill] at (9,4) {};
\node[draw,circle,inner sep=1.5pt, fill] at (9,6) {};
\node[draw,circle,inner sep=1.5pt, fill] at (11,2) {};
\node[draw,circle,inner sep=1.5pt, fill] at (11,4) {};
\node[draw,circle,inner sep=1.5pt, fill] at (11,6) {};
\node[blue,draw,circle,inner sep=2pt, fill] at (8,8) {};
\node[blue,draw,circle,inner sep=2pt, fill] at (10,8) {};
\node[blue,draw,circle,inner sep=2pt, fill] at (12,8) {};
\node[magenta,draw,circle,inner sep=2pt, fill] at (14,7) {};
\node[violet,draw,circle,inner sep=2pt, fill] at (14,3) {};
\node[green,draw,circle,inner sep=2pt, fill] at (14,1) {};
\node[orange,draw,circle,inner sep=2pt, fill] at (14,5) {};
\node[right] at (14.3,7) {1};
\node[right] at (14.3,5) {-1};
\node[right] at (14.3,3) {2};
\node[right] at (14.3,1) {0};
\draw[gray] (7,6)--(14,7)--(9,6); 
\draw[gray] (11,6)--(14,7)--(7,4);
\draw[gray] (11,4)--(14,7);
\draw[gray] (7,2)--(14,3);
\draw[gray] (9,4)--(14,5)--(11,2);
\draw[gray] (9,2)--(14,1);
\draw[gray] (7,4)--(8,8)--(7,6);
\draw[gray] (7,2)--(8,8);
\draw[gray] (9,4)--(10,8)--(9,6);
\draw[gray] (9,2)--(10,8);
\draw[gray] (11,4)--(12,8)--(11,6);
\draw[gray] (11,2)--(12,8);
\draw[gray] (9,2)--(4,3)--(7,2);
\draw[gray] (11,2)--(4,3);
\draw[gray] (9,4)--(4,5)--(7,4);
\draw[gray] (11,4)--(4,5);
\draw[gray] (9,6)--(4,7)--(7,6);
\draw[gray] (11,6)--(4,7);
\end{tikzpicture}
\end{center}
\caption{The graph corresponding to the character table of the symmetric group on 3 symbols.}
\end{figure}
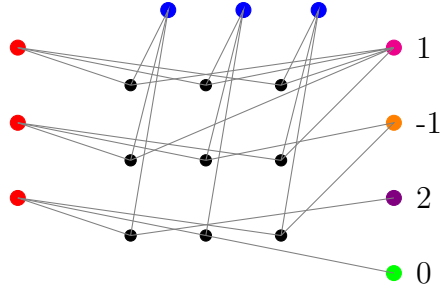

We can see how permuting the rows and columns of $\CT(S)$ does not change the graph $\mathcal{S}$. Thus the graph $\CT(S)$ captures the equivalence seen in equation \ref{eqn 1}.

\subsection{Canonically label the graph.}

To determine a canonical labeling of the graph we used McKay's nauty algorithm \cite{nauty} as implemented in the {\sc GAP} \cite{GAP} package {\sc Digraphs} \cite{Digraphs} (and by association the package {\sc Grape} \cite{Grape}). The algorithm works by manipulating the symmetry of the table. 

\subsection{Hash the canonical labeling of $\mathcal{G}$.}
We cast the canonically labeled graph as a string and then output the \texttt{md5sum} \cite{MD5} of the string.

\subsection{Sort the list.}

We sorted the list groups, designated by their indexes in the SmallGroup database, and collected those with identical hash and number of conjugacy classes into potential equivalence classes.  

\subsection{Check for hash collisions.}\label{sec collisions}

To check for collisions, we used the function \texttt{IsIsomorphic(CT(G),CT(H))} in {\sc Magma}  \cite{Magma} for groups $G$ and $H$ in the same potential equivalence class. These calculations could have been quadratic in the size of the potential equivalence class. However, because we discovered no collisions from the hash function, meaning that the hash of the canonical labeling of $\mathcal{G}$ determined the equitabular class of $G$, we only had to compare each member of the partition against a single fixed group to verify that our partition was in fact an equitabular class. 

We chose to run this portion of the computation in {\sc Magma} to give a complementary check to the partition coming from {\sc GAP}. 

\subsection{Comparison with the general approach.}\label{sec compare}
In general, the historical approach to the problem of determining all of the equitabular classes for a given order $n$ was to first compute a number of invariants on all of the groups of order $n$. These invariants were information about a group $G$ that is preserved by $\CT(G)$, meaning that the value of the invariant could be read from a table representing $\CT(G)$. Among these invariants are the order of $G$, the number of conjugacy classes of $G$, and the size of the conjugacy classes. In general it is not known what list of invariants is sufficient to identify the equitabular classes of groups of a given order. 

Ultimately, one is forced to compare the character tables of groups that sometimes do not have the same character table. For orders with a small number of groups this is not computationally expensive. In our approach to handle orders with a large number of groups, e.g., 1536, we ended up only comparing groups with the same character table---the ``golden'' invariant we stumbled across was the hash of the canonical labeling of the table. By utilizing high-throughput computing, computing this value for all groups of a given order was feasible.

We wish to point out that high-throughput computing involves distributing the problem to multiple machines each of which work independently of each other. In contrast the high-performance-computing model runs jobs on a tightly coupled cluster of machines which have extensive CPU and RAM to perform calculations. We crafted our algorithm to work on available high-throughout-computing resources. 

In implementing our algorithm each character table was hashed independently to produce a small output file which allowed for easy file transfer and low requirements on the CPU or RAM. This enabled us to access a majority of nodes on in our computing cluser. Our algorithm was successful because: \begin{itemize} \item it had a small output for each job; \item it utilized open source software; \item it had low system requirements. \end{itemize} 

\section{Examples and observations.}\label{sec observations}
In this section, we provide some examples and observations of how one could query the database. 

\subsection{Derived Length}\label{sec derived length}
Using the character table of a finite group, one can determine if the group is solvable by constructing the lattice of normal subgroups of $G$ and looking for a chain of normal subgroups such that the index of each subgroup in the next is a prime power. Moreover, the conjugacy classes contained in the derived subgroup can be readily identified from the character table, as can the commutators. However, the derived length of the group is not observable from the character table by itself. 

Mattarei first observed that the derived lengths cannot be determined from the character table of a group; meaning that there are groups $G$ and $H$ with the same character table, but with different derived lengths \cite{Mattarei1}. In later work he produced $p$-groups $G$ and $H$ with the same character table, but different derived lengths \cite{Mattarei2}. For his $p$-group examples, Mattarei required that $p\geq 5$ and the groups themselves have order $p^{11}$. Using our database, we were able to identify the following examples of groups of order $512=2^9$ that have the same character table, but different derived lengths; moreover, we can easily verify that they are the smallest possible examples in terms of order. 

\begin{example}
This pseudocode could be used to see that there are no groups with order less than 512 that have the same character table and different derived lengths.
\begin{verbatim}
for n  in 1 to 511 do
      for E in equitabular_classes(n) do 
            for G and H in E do
                  if DerivedLength(G) ne DerivedLength(H) then 
                        print Index(G), Index(H).
                  end if;
            end for;
      end for;
end for;
\end{verbatim}
However by calculating the derived lengths for all equitabular groups of order 512, we get the following table. 
\begin{table} [h!] 
\begin{center}
\begin{tabular} {l  l l }
 &Derived Length 2 & Derived Length 3 \\
\hline
Class 1 &1637, 1638, 1639, 1640 & 1615, 1616, 1617, 1618, \\ &&1619, 1620, 1621, 1622\\\hline
Class 2 &46947, 46948, 46949, 46950, & 46929, 46930, 46931, 46932, \\
&46951, 46952, 46953, 46954 &46933, 46934, 47002, 47003, \\&&47004 \\\hline
Class 3 &59243, 59244, 59245, 59246, &59213, 59214, 59215, 59216, \\
&59247, 59248 & 59217, 59218, 59219, 59220 \\\hline
Class 4 &59930, 59931, 59932, 59933, & 59906, 59907, 59908, 59909, \\
&59934, 59935, 59936, 59937 & 59910, 59911, 59912, 59913\\
\end{tabular}
\caption{Indices of groups of order 512 whose equitabular class contains groups with different derived lengths.}
\end{center}
\end{table}

\end{example}

We mention that it is unknown if there are equitabular groups whose derived lengths differ by more than one. 

\subsection{Character Theoretic Words}

Another area of recent interest in group theory is whether the image of certain word maps can be identified from the character table of a group $G$.  It is well-known that the conjugacy classes of a group $G$ whose elements occur as commutators can be identified from the character table of $G$, i.e., looking at the character values over $g$ we can tell if there is some $x,y\in G$ with $[x,y]=x^{-1}y^{-1} x y = g$. Explicitly, we have the following lemma, which is an exercise in \cite[Exercise 3.10 ]{Isaacs_CT} and a lemma in \cite[Lemma 2.6.4]{LP}.

\begin{lemma}
Let $G$ be a finite group and let $g\in G$. There is some $x,y\in G$ with $g=[x,y]$ if and only if 
\[
\sum_{\chi \in \mathrm{Irr}(G)} \frac{\chi(g)}{\chi(1)} \ne 0.
\]
\end{lemma}

The character table can also be used to determine the number of ways $g$ occurs as a commutator. The recently proven Ore conjecture asked whether every element of a finite nonabelian simple group $G$ occurred as a commutator. The proof of the Ore conjecture by Liebeck, O'Brien, Shalev, and Tiep utilized the character theoretic nature of the word $w=x^{-1} y^{-1} xy$ \cite{LOST}. There has been some interest over finite nonabelian simple groups of the probability that $g$ occurs as a commutator \cite{GS}; this probability is also determined by the character table of $G$. 

Besides $[x,y]$, there are results known for other words, for example $w=x^2y^2$. The word $w$ is also character theoretic in that:
\begin{lemma}\cite[Lemma 2.2]{LOST2}
Let $G$ be a finite group and $g\in G$. The number of ways $g$ occurs as a product of two squares is 
\[
|G| \cdot \sum_{\substack{\chi \in \mathrm{Irr}(G)\\ \chi \mathrm{\,\,  real}}} \frac{\chi(g)}{\chi(1)}.
\]
\end{lemma}
The authors who proved the Ore conjecture also showed that every element of a finite nonabelian simple group is a product of two squares \cite{LOST2}. 

We will say that the image of a word $w \in \textbf{F}_n$ is the set of all evaluations of $w$ in $G$. A natural question to ask is whether for every word $w$, there is some way to deduce from the character table of a group $G$ whether or not an element $g\in G$ occurs in the image of $w$. As seen above when $w=[x,y]$ or $w=x^2y^2$ this is the case. However, when $w=x^p$ for $p$ an odd prime, we note that the extraspecial groups of order $p^3$ share a character table. Because one of the extraspecial groups of order $p^3$ has exponent $p$ and one has exponent $p^2$, the image of $w$ cannot be extracted from the character table by itself. 

What about the word $w=x^2$?

\begin{example}
We can ask the character table database to look for groups $G$ and $H$ such $G$ and $H$ have the same character table, but the number of squares in $G$ is different than the number in $H$. 
\begin{verbatim}
for n in 1 to 1000 do
      for E in equitabular_classes(n) do
            for G and H in E do
                  if #Squares in G ne #Squares in H then 
                        print n, Index(G), Index(H).
                  end if;
            end for;
      end for;
end for;		
\end{verbatim}
An algorithm based on the pseudocode above returns that when $G=\mathrm{SmallGroup(64,100)}$ and $H=\mathrm{SmallGroup(64,98)}$, then $G$ and $H$ share a character table and $|\{x^2: x \in G\}|=6$ and $|\{x^2: x \in H\}|=5$. 
\end{example}

From the above example, we see that $\CT(G)$ determines the group generated by the word $w=x^2$, but not the image of the word itself. 

\section{Availabiltiy.}\label{sec availability}
The database can be accessed at \url{https://osf.io/7zjh9/}. In addition we have included some tools to query the database and convert the data to a format usable in either {\sc Magma} or {\sc GAP}.

\section{Acknowledgements.}\label{sec acknowledgements}
We want to thank Gerhard Hiss for helping us obtain a copy of Skrzipczyk's thesis. 

We also want to thank James Mitchell and Wilf Wilson for helping with the {\sc Digraphs} package.

This research was performed using the compute resources and assistance of the UW-Madison Center For High Throughput Computing (CHTC) in the Department of Computer Sciences. The CHTC is supported by UW-Madison, the Advanced Computing Initiative, the Wisconsin Alumni Research Foundation, the Wisconsin Institutes for Discovery, and the National Science Foundation, and is an active member of the Open Science Grid, which is supported by the National Science Foundation and the U.S. Department of Energy's Office of Science. 

This research was done using resources provided by the Open Science Grid \cite{OSG1,OSG2}, which is supported by the National Science Foundation award 1148698, and the U.S. Department of Energy's Office of Science.

This material is based upon work done while the first author was supported by the National Science Foundation under Grant No. DMS-1502553. The first author also acknowledges that this material is based upon work supported by the National Science Foundation Graduate Research Fellowship Program under Grant No. DGE-1256529.

\bibliographystyle{alpha}
\bibliography{ref}

\end{document}